\documentstyle[11pt]{article}

\author{Charles H. Walter\thanks{Supported in part by NSA research grant
MDA904-92-H-3009.}\\
Laboratoire de Math\'ematiques\\
Universit\'e de Nice\\
F-06108 Nice Cedex 02 FRANCE}
\title{On the Harder-Narasimhan Filtration for Coherent Sheaves on
${\bf P}^2$: I}
\date{ }

\newtheorem{theorem}[subsection]{Theorem}
\newtheorem{lemma}[subsection]{Lemma}

\newtheorem{defn}[subsection]{\sl Definition}
\newenvironment{definition}{\begin{defn}\rm}{\end{defn}}

\input amssym.def
\newsymbol\subsetneqq 2324

\def\qed{\hfill $\Box$}
\long\def\TeXButton#1#2{#2}
\def\proof{\paragraph{Proof. }}

\catcode`\@=11

\def\@begintheorem#1#2{\sl \trivlist
   \item[\hskip \labelsep{\bf #1\ #2\thmcounterend}]}
\def\@opargbegintheorem#1#2#3{\sl \trivlist
      \item[\hskip \labelsep{\bf #1\ #2\ (#3)\thmcounterend}]}
\def\thmcounterend{.}

\def\section{\@startsection{section}{1}{\z@}{-3.25ex plus
 -1ex minus -.2ex}{1.5ex plus .2ex}{\large\bf}}
\def\subsection{\@startsection
 {subsection}{2}{\z@}{3.25ex plus 1ex minus .2ex}{-0.5em}{\normalsize\sl}}
\def\subsubsection{\@startsection
 {subsubsection}{3}{\z@}{3.25ex plus 1ex minus .2ex}{-0.5em}{\normalsize\sl}}
\def\paragraph{\@startsection
 {paragraph}{3}{\z@}{2ex plus 0.6ex minus .2ex}{-0.5em}{\normalsize\sl}}
\def\subparagraph{\@startsection
 {subparagraph}{3}{\parindent}{2ex plus 0.6ex minus
.2ex}{-1pt}{\normalsize\sl}}

\@addtoreset{equation}{subsection}

\def\RIfM@{\relax\ifmmode}

\newif\iffirstchoice@
\firstchoice@true
\def\textfonti{\the\textfont\@ne}
\def\textfontii{\the\textfont\tw@}
\def\text{\RIfM@\expandafter\text@\else\expandafter\text@@\fi}
\def\text@@#1{\leavevmode\hbox{#1}}
\def\text@#1{\mathchoice
 {\hbox{\everymath{\displaystyle}\def\textfonti{\the\textfont\@ne}%
  \def\textfontii{\the\textfont\tw@}\textdef@@ T#1}}
 {\hbox{\firstchoice@false
  \everymath{\textstyle}\def\textfonti{\the\textfont\@ne}%
  \def\textfontii{\the\textfont\tw@}\textdef@@ T#1}}
 {\hbox{\firstchoice@false
  \everymath{\scriptstyle}\def\textfonti{\the\scriptfont\@ne}%
  \def\textfontii{\the\scriptfont\tw@}\textdef@@ S\rm#1}}
 {\hbox{\firstchoice@false
  \everymath{\scriptscriptstyle}\def\textfonti
  {\the\scriptscriptfont\@ne}%
  \def\textfontii{\the\scriptscriptfont\tw@}\textdef@@ s\rm#1}}}
\def\textdef@@#1{\textdef@#1\rm\textdef@#1\bf\textdef@#1\sl\textdef@#1\it}
\def\DN@{\def\next@}
\def\eat@#1{}
\def\textdef@#1#2{%
 \DN@{\csname\expandafter\eat@\string#2fam\endcsname}%
 \if S#1\edef#2{\the\scriptfont\next@\relax}%
 \else\if s#1\edef#2{\the\scriptscriptfont\next@\relax}%
 \else\edef#2{\the\textfont\next@\relax}\fi\fi}

\def\tsum{\mathop{\textstyle \sum }}
\def\binom#1#2{{#1 \choose #2}}

\catcode`\@=12

\begin{document}

\maketitle
\begin{abstract}
\noindent Let ${\cal E}$ be a torsion-free sheaf on ${\bf P}^2$. We give an
effective method which uses the Hilbert function of ${\cal E}$ to construct
a weak version of the Harder-Narasimhan filtration of a torsion-free sheaf
on ${\bf P}^2$ subject only to the condition that ${\cal E}$ be sufficiently
general among sheaves with that Hilbert function. This algorithm uses on a
generalization of Davis' decomposition lemma to higher rank. \bigskip\
\end{abstract}

\TeXButton{subsection1}{\refstepcounter{subsection}}Consider the following
problem. Let ${\cal E}$ be an explicit torsion-free sheaf on ${\bf P}^2$
given by a presentation
\begin{equation}
\label{pres}0\rightarrow \bigoplus_{n\in {\bf Z}}{\cal O}_{{\bf P}%
^2}(-n)^{b(n)}\stackrel{\phi }{\rightarrow }\bigoplus_{n\in {\bf Z}}{\cal O}%
_{{\bf P}^2}(-n)^{a(n)}\rightarrow {\cal E}\rightarrow 0.
\end{equation}
How does one go about effectively computing the Harder-Narasimhan filtration
of ${\cal E}$, i.e.\ the unique filtration%
$$
0=F_0({\cal E})\subset F_1({\cal E})\subset \cdots \subset F_s({\cal E})=%
{\cal E}
$$
such that the graded pieces ${\rm gr}_i({\cal E}):=F_i({\cal E})/F_{i-1}(%
{\cal E})$ are semistable in the sense of Gieseker-Maruyama and their
reduced Hilbert polynomials $P_i(n)=\chi ({\rm gr}_i({\cal E})(n))/{\rm rk}(%
{\rm gr}_i({\cal E}))$ satisfy $P_1(n)>P_2(n)>\cdots >P_s(n)$ for $n\gg 0?$

In this paper and its planned sequel we consider the problem under the
simplifying assumption that the matrix $\phi $ of homogeneous polynomials is
general, i.e.\ that ${\cal E}$ is general among torsion-free sheaves with
the same Hilbert function as ${\cal E}$. Our solution to the problem then
divides into two parts. In this first part we construct a filtration of $%
{\cal E}$ of the type
\begin{equation}
\label{genfilt}0\subset {\cal E}_{\leq \tau _1}\subset \cdots \subset {\cal E%
}_{\leq \tau _s}\subset {\cal E}
\end{equation}
where ${\cal E}_{\leq n}$ denotes the subsheaf of ${\cal E}$ which is the
image of the natural evaluation map $H^0({\cal E}(n))\otimes {\cal O}_{{\bf P%
}^2}(-n)\rightarrow {\cal E}$. We give an algorithm for picking the $\tau _i$
so that the filtration approximates the true Harder-Narasimhan filtration
but groups together all pieces of the Harder-Narasimhan filtration with
slopes between two consecutive integers. The associated graded sheaves ${\rm %
gr}_i({\cal E}):={\cal E}_{\leq \tau _i}/{\cal E}_{\leq \tau _{i-1}}$ are
not always semistable, but they do share a number of properties with
semistable sheaves which suffice for a number of applications. For instance
they are of {\em rigid splitting type}, i.e.\ their restrictions to a
general line $L$ of ${\bf P}^2$ are of the form ${\cal O}_L(n_i)^{\alpha
_i}\oplus {\cal O}_L(n_i+1)^{\beta _i}$ for some $n_i$, $\alpha _i$ and $%
\beta _i$. They also satisfy ${\rm Hom}({\rm gr}_i({\cal E}),{\rm gr}_i(%
{\cal E})(-1))=0$.

We call our filtration the Weak Harder-Narasimhan (or WHN) filtration of $%
{\cal E}$. It is fine enough to give graded pieces of rigid splitting type
but is otherwise deliberately as coarse as possible in order to keep the
algorithm for picking the $\tau _i$ as simple as possible (and also because
further refinement can actually be counterproductive in such applications as
the classification of irreducible components of the moduli stack of
torsion-free sheaves on ${\bf P}^2$). Thus in some cases the WHN filtration
may not even be the finest filtration of ${\cal E}$ by subsheaves of the
form ${\cal E}_{\leq n}$ which is compatible with the Harder-Narasimhan
filtration, although such a refinement could certainly be computed by the
methods of this paper by adding an extra step to the end of the algorithm of
paragraph (\ref{effective}). The true Harder-Narasimhan filtration is not
always given by subsheaves of the form ${\cal E}_{\leq n}$ and can therefore
be much harder to compute. For example the sheaf ${\cal E=O}_{{\bf P}%
^2}\oplus \Omega _{{\bf P}^2}(2)$ has Harder-Narasimhan filtration $0\subset
\Omega _{{\bf P}^2}(2)\subset {\cal E}$ but is unfilterable by subsheaves of
the form ${\cal E}_{\leq n}$ since ${\cal E}_{\leq n}=0$ for $n\leq -1$ and $%
{\cal E}_{\leq n}={\cal E}$ for $n\geq 0$.

In the planned part II we will show how to refine the WHN filtration of a
sufficiently general sheaf to the true Harder-Narasimhan filtration using
exceptional objects and mutations.

The precise formulation of the WHN filtration requires a certain number of
numerical definitions. We consider a general sheaf ${\cal E}$ with a
presentation of the form (\ref{pres}) for given functions $a(n)$ and $b(n)$
of finite support. We define $r(n)$ and $h(n)$ as the first and second
integrals of $a(n)-b(n)$, i.e.%
\begin{eqnarray}
r(n) & := & \sum_{m\leq n}\left\{ a(m)-b(m)\right\} ,
\label{r(n)} \\ h(n) & := & \sum_{m\leq n}r(m)=\sum_{m\leq n}(n-m+1)\left\{
a(m)-b(m)\right\} .\label{h(n)}
\end{eqnarray}The function $h$, $r$, and $a-b$ are respectively the first,
second, and third differences of the Hilbert function of ${\cal E}$ defined
by $n\mapsto h^0({\cal E}(n))$. We will assume that the $a(n)$ and $b(n)$
are such that $r(n)\geq 0$ for all $n$. The general $\phi {:}~{\cal %
\bigoplus O}(-n)^{b(n)}\rightarrow \bigoplus {\cal O}(-n)^{a(n)}$ is
injective if and only if this is the case (see \cite{Ch} or Theorem \ref
{Chang} below). Depth considerations show that the cokernel ${\cal E}$ of
such an injective $\phi $ will have no subsheaves supported at isolated
points, but ${\cal E}$ is permitted to have torsion supported along a curve.

We now define further auxiliary functions by%
\begin{eqnarray}
\tilde h(n) & := & \max \, \{h(m)+(n-m)r(m) \mid m\geq n \}, \label{htilde(n)}
\\
t(n) & := & \max \,\{m\geq n\mid \tilde
h(n)=h(m)-(m-n)r(m)\}\in {\bf Z}\cup \{+\infty \}.\label{t(n)}
\end{eqnarray}We will show in Lemma \ref{combin} that if $a(n)$ and $b(n)$
are such that $r(n)\geq 0$ for all $n$, then the function $t$ is
nondecreasing and takes only finitely many values $\tau _0<\tau _1<\cdots
<\tau _s<\tau _{s+1}=+\infty $. These $\tau _i$ may be effectively computed
by an algorithm we will give in paragraph (\ref{effective}). We set also $%
\tau _{-1}=-\infty $. Then we define the {\em WHN filtration} of ${\cal E}$
as the filtration
\begin{equation}
\label{wfilt}0={\cal E}_{\leq \tau _{-1}}\subset {\cal E}_{\leq \tau
_0}\subset \cdots \subset {\cal E}_{\leq \tau _s}\subset {\cal E}_{\leq \tau
_{s+1}}={\cal E}
\end{equation}
with graded pieces ${\rm gr}_i({\cal E}):={\cal E}_{\leq \tau _i}/{\cal E}%
_{\leq \tau _{i-1}}$ for $0\leq i\leq s+1$. Our main result is:

\begin{theorem}
\label{main}Let $a$, $b{:}~{\bf Z}\rightarrow {\bf Z}_{\geq 0}$ be functions
of finite support such that the function $r(n)$ of (\ref{r(n)}) is
nonnegative. Let ${\cal E}$ be the cokernel of an injection $\phi %
{:}~\bigoplus_n{\cal O}_{{\bf P}^2}(-n)^{b(n)}\rightarrow \bigoplus_n%
{\cal O}_{ {\bf P}^2}(-n)^{a(n)}$. If $\phi $ is sufficiently general, then
the WHN filtration of ${\cal E}$ defined in (\ref{wfilt}) has the following
properties:

(i)\quad For all $0\leq i\leq s+1$ the sheaf ${\rm gr}_i({\cal E})$ has
resolution%
$$
0\rightarrow \bigoplus_{\tau _{i-1}<n\leq \tau _i}{\cal O}_{{\bf P}%
^2}(-n)^{b(n)}\rightarrow \bigoplus_{\tau _{i-1}<n\leq \tau _i}{\cal O}_{%
{\bf P}^2}(-n)^{a(n)}\rightarrow {\rm gr}_i({\cal E})\rightarrow 0.
$$

(ii)\quad The subsheaf ${\cal E}_{\leq \tau _0}={\rm gr}_0({\cal E})$ is the
torsion subsheaf of ${\cal E}$.

(iii)\quad For $1\leq i\leq s+1$ the sheaf ${\rm gr}_i({\cal E})$ is
torsion-free and of rigid splitting type, i.e.\ if $L$ is a general line of $%
{\bf P}^2$, then ${\rm gr}_i({\cal E}){\mid }_L\cong {\cal O}_L(-\nu
_i)^{\beta _i}\oplus {\cal O}_L(-\nu _i-1)^{\rho _i-\beta _i}$ for some
integers $\nu _i$, $\beta _i$ and $\rho _i$. Moreover $\nu _1<\nu _2<\cdots
<\nu _{s+1}$ and ${\cal E}{\mid }_L\cong \bigoplus_{i=0}^{s+1}{\rm gr}_i(%
{\cal E}){\mid }_L$.

(iv)\quad For $i\leq j$ we have ${\rm Hom}({\rm gr}_i({\cal E}),{\rm gr}_j(%
{\cal E})(-1))=0$.

(v)\quad The Harder-Narasimhan filtration of ${\cal E}/{\cal E}_{\leq \tau
_0}$ for Gieseker-Maruyama stability is a refinement of the filtration (\ref
{wfilt}) of ${\cal E}/{\cal E}_{\leq \tau _0}$. Indeed ${\rm gr}_i({\cal E})$
collects all pieces of the Harder-Narasimhan filtration with slopes $\mu $
satisfying $-\nu _i-1<\mu <-\nu _i$ as well as some of those of slopes $-\nu
_i-1$ and $-\nu _i$.
\end{theorem}

The outline of the paper is as follows. In the first section we prove a
number of numerical lemmas leading to a method of filtering the Hilbert
function of ${\cal E}$. The key definition is that the Hilbert function of $%
{\cal E}$ (or its differences $h(n)$, $r(n)$ or $a(n)-b(n)$ as defined
above) is {\em filterable} at $m$ if the function $r(n)$ of (\ref{r(n)})
satisfies $r(n)\geq r(m)$ for all $n\geq m$. The $a(n)$ and $b(n)$ then
split into
$$
a_m^{{\rm sub}}(n):=\left\{
\begin{array}{ll}
a(n) & \text{if }n\leq m, \\ 0 & \text{if }n\geq m,
\end{array}
\right.
$$
and analogous functions $a_m^{{\rm quot}}(n)$, $b_m^{{\rm sub}}(n)$ and $%
b_m^{{\rm quot}}(n)$. The $r_m^{{\rm sub}}$, $r_m^{{\rm quot}}$, $h_m^{{\rm %
sub}}$ and $h_m^{{\rm quot}}$ are defined by integrating. If a Hilbert
function is filterable at several integers $m_i$, it may be split into
several graded pieces this way. The lemmas of the section show that the
Hilbert function of ${\cal E}$ is filterable at the $\tau _i$ and that its
graded pieces satisfy conditions analogous to the conditions of parts
(ii)-(v) of Theorem \ref{main}.

In the second section we show that such filtrations of Hilbert functions
correspond to filtrations of ${\cal E}$ by subsheaves of the form ${\cal E}%
_{\leq m}$ if ${\cal E}$ is sufficiently general among coherent sheaves with
the same Hilbert function. The key lemma is the following which may be
regarded as a generalization of Davis' decomposition lemma \cite{D} to
higher rank.

\TeXButton{Davis.lemma}
{\begin{trivlist}
\item [\hskip \labelsep{\bf Lemma \ref{gener}.}]{\sl
Suppose that ${\cal E}$ is a coherent sheaf on ${\bf P}^2$ without
zero-dimensional associated points such that the Hilbert function of ${\cal E}$
is filterable at an integer $m$. Write ${\cal E}$ as the cokernel of an
injection $\phi {:}\ \bigoplus_n{\cal O}_{{\bf P}^2}(-n)^{b(n)}\rightarrow
\bigoplus_n{\cal O}_{{\bf P}^2}(-n)^{a(n)}$. If the matrix $\phi $ is
sufficiently general, then ${\cal E}_{\leq m}$ and ${\cal E}/{\cal E}_{\leq
m}$ have resolutions}\begin{eqnarray*}
& 0\rightarrow \bigoplus_{n\leq m}{\cal O}_{{\bf P}^2}(-n)^{b(n)}\rightarrow
\bigoplus_{n\leq m}{\cal O}_{{\bf P}^2}(-n)^{a(n)}\rightarrow {\cal E}_{\leq
m}\rightarrow 0 , &  \\
& 0\rightarrow \bigoplus_{n>m}{\cal O}_{{\bf P}^2}(-n)^{b(n)}\rightarrow
\bigoplus_{n>m}{\cal O}_{{\bf P}^2}(-n)^{a(n)}\rightarrow {\cal
E}_{>m}\rightarrow 0 . &
\end{eqnarray*}
\end{trivlist}
}

The rest of the section is devoted to showing that Theorem \ref{main}
follows from this lemma and from the numerical lemmas proved in the first
section.

This paper was written in the context of the group on vector bundles on
surfaces of Europroj. The author would like to thank A.\ Hirschowitz for
some useful conversations.

\section{Filtering Hilbert Functions}

This section contains the purely combinatorial part of the proof of the
Theorem \ref{main}. It consists of a number of numerical lemmas on Hilbert
functions of coherent sheaves on ${\bf P}^2$. We begin by fixing some
terminology. We use the notation $(x)_{+}=\max (x,0)$.

Our fundamental invariant is the difference $a(n)-b(n)$ between the
functions of (\ref{pres}). We assume that $a(n)-b(n)$ is an integer for all $%
n$ and vanishes for all but finitely many $n$ and that the associated
function $r(n)=\sum_{m\leq n}\left\{ a(n)-b(n)\right\} \geq 0$ for all $n$.
We call the associated function $h(n)=\sum_{m\leq n}r(n)$ of (\ref{h(n)})
the {\em FDH function} (or first difference of a Hilbert function). It is
the FDH functions which will play the major role in our computations. An
intrinsic definition is:

\begin{definition}
An FDH function is a function $h{:}~{\bf Z}\rightarrow {\bf Z}_{\geq 0}$
such that $r(n)=\Delta h(n)\geq 0$ for all $n$, $h(n)=0$ for $n\ll 0$, and $%
h(n)$ is linear of the form $\rho n+\sigma $ for $n\gg 0$. We call $\rho $
the {\em rank} of $h$, $\sigma -\rho $ its {\em degree}, and $\sum_{n\in
{\bf Z}}\left\{ (\rho n+\sigma )_{+}-h(n)\right\} $ its {\em deficiency}.

An FDH function is {\em torsion} if its rank is $0$. An FDH function $h$ is
{\em torsion-free} if $r(m)\geq 1$ implies that $r(n)\geq 1$ for all $n\geq m
$. An FDH function $h$ is {\em locally free} if $h$ is torsion-free and
additionally $r(m)\geq 2$ implies that $r(n)\ge 2$ for all $n\geq m$. (This
terminology will be justified by Theorem \ref{Chang}.)
\end{definition}

These functions have the following basic properties:

\begin{lemma}
\label{combin}Suppose $h$ is an FDH function of rank $\rho $ and degree $%
\sigma -\rho $. Let $r=\Delta h$, and let $\tilde h$ and $t$ be as in (\ref
{htilde(n)}) and (\ref{t(n)}). Then

(i)\quad For all $n$ one has $0\leq h(n)\leq \tilde h(n)$ and $t(n)>n$.

(ii)\quad If $m>t(n)$, then $r(m)>r(t(n))$.

(iii)\quad The function $t$ is nondecreasing and takes only finitely many
distinct values. If we write these as $\tau _0<\tau _1<\cdots <\tau _s<\tau
_{s+1}=\infty $, then $0=r(\tau _0)<r(\tau _1)<\cdots <r(\tau _s)<\rho $.

(iv)\quad Let $\nu _i=\min \{n\mid t(n)=\tau _i\}$. If $\nu _i\leq n<\nu
_{i+1}$, then $t(n)=\tau _i$ and $\tilde h(n)=nr(\tau _i)+\left( h(\tau
_i)-\tau _ir(\tau _i)\right) $. Moreover $n<\nu _i$ if and only if
\begin{equation}
\label{nu(i)}nr(\tau _{i-1})+\left( h(\tau _{i-1})-\tau _{i-1}r(\tau
_{i-1})\right) >nr(\tau _i)+\left( h(\tau _i)-\tau _ir(\tau _i)\right)
\end{equation}

(v)\quad If $\nu _i\leq n\leq m\leq \tau _i$, then $h(\tau _i)+(n-\tau
_i)r(\tau _i)\geq h(m)+(n-m)r(m)$.

(vi)\quad If $n<\nu _1$, then $t(n)=\tau _0=\max \{n\mid r(n)=0\}$ and $%
\tilde h(n)=h(\tau _0)$. In particular, if $h$ is torsion-free then $\tilde
h(n)=0$ and $t(n)=\max \{m\mid h(m)=0\}$ for $n<\nu _1$.

(vii)\quad For $n\geq \nu _{s+1}$ one has $\tilde h(n)=\rho n+\sigma $ and $%
t(n)=\tau _{s+1}=+\infty $.
\end{lemma}

\TeXButton{Proof}{\proof}(i) From the definitions we see that $%
h(n)=h(n+1)-r(n+1)\leq \tilde h(n)$ and that this implies that $t(n)\geq n+1$%
{}.

(ii) We go by induction on $m$. Thus we assume that $r(i)>r(t(n))$ for $%
t(n)<i<m$, and we will show that $r(m)>r(t(n))$ as well. But the definitions
of $\tilde h(n)$ and $t(n)$ yield immediately%
$$
h(t(n))+(n-t(n))r(t(n))=\tilde h(n)>h(m)-(m-n)r(m).
$$
Hence%
\begin{eqnarray*}
(m-n)r(m) & > & (t(n)-n)r(t(n))+h(m)-h(t(n)) \\
& = & (t(n)-n)r(t(n))+\tsum_{i=t(n)+1}^mr(i) \\
& \geq  & (t(n)-n)r(t(n))+(m-t(n)-1)r(t(n))+r(m)
\end{eqnarray*}whence (ii).

(iii) Since $t(n-1)\geq n$ by (i), we see from the definitions that
$$
\tilde h(n)\geq h(t(n-1))+(n-t(n-1))r(t(n-1))=\tilde h(n-1)+r(t(n-1)).
$$
Thus $r(t(n-1))\leq \Delta \tilde h(n)$. Similarly%
$$
\tilde h(n-1)\geq h(t(n))+(n-1-t(n))r(t(n))=\tilde h(n)-r(t(n))
$$
and $\Delta \tilde h(n)\leq r(t(n))$. Hence $r\circ t$ is nondecreasing.
Because of (ii) this implies that $t$ is nondecreasing. The function $t$ can
only take finitely many values since by (ii) $r$ takes a different value at
each value of $t$, and the values of $r$ are bounded since $r(n)$ is
constant for $n\ll 0$ and for $n\gg 0$. If $\tau _i=$$t(n)<+\infty $, then
for $m\gg \tau _i$ (ii) yields $\rho =r(m)>r(\tau _i)$. Finally for $n\ll 0$
one has $0\leq \tilde h(n)=h(\tau _0)+(n-\tau _0)r(\tau _0)$ which implies
that $r(\tau _0)\leq 0$. But since $h$ is an FDH function, $r(\tau _0)\geq 0$%
. So $r(\tau _0)=0$. This finishes (iii).

(iv) If $\nu _i\leq n<\nu _{i+1}$, then $t(n)=\tau _i$. According to the
definitions, this implies
$$
\tilde h(n)=h(\tau _i)-(\tau _i-n)r(\tau _i)=nr(\tau _i)+\left( h(\tau
_i)-\tau _ir(\tau _i)\right) .
$$
As for the inequality (\ref{nu(i)}), because both its sides are linear and
the slope on the left side is less than that on the right side, it is enough
to show that the inequality holds for $n=\nu _i-1$ but fails for $n=\nu _i$.
But because of the definition of $t$, this follows immediately from $t(\nu
_i-1)=\tau _{i-1}<$$\tau _i=t(\nu _i)$.

(v) The proof is divided into several cases. First if $n<\nu _{i+1}$, then
by (iv) the inequality becomes $\tilde h(n)\geq h(m)-(m-n)r(m)$ which
follows from the definition of $\tilde h(n)$. If $n\ge \nu _{i+1}$ but $%
r(m)\leq r(\tau _i)$, then the inequality follows from the case $n=\nu _i$ by%
\begin{eqnarray*}
h(\tau _i)+(n-\tau _i)r(\tau _i) & = & \tilde h(\nu _i)+(n-\nu _i)r(\tau _i) \\
  & \geq  & h(m)+(\nu _i-m)r(m)+(n-\nu _i)r(\tau _i) \\
  & \geq  & h(m)+(n-m)r(m).
\end{eqnarray*}Finally if $n\ge \nu _{i+1}$ but $r(m)>r(\tau _i)$, then let $%
m^{\prime }:=\min \{M>m\mid r(M)\leq r(\tau _i)\}\leq \tau _i$. Then using
the previous case applied with $m^{\prime }-1$ substituted for $n$ and $%
m^{\prime }$ substituted for $m$ we see that%
\begin{eqnarray*}
h(m)+(n-m)r(m) & = & \left[ h(m^{\prime })-r(m^{\prime })\right] -
\sum_{i=m+1}^{m^{\prime }-1}r(i)-(m-n)r(m) \\
  & \leq  & \left[ h(\tau _i)+(m^{\prime }-1-\tau _i)r(\tau _i)\right] -
(m^{\prime }-1-n)r(\tau _i) \\
  & = & h(\tau _i)+(n-\tau _i)r(\tau _i) .
\end{eqnarray*}

(vi) If $n<\nu _1$, then $t(n)=\tau _0$ by (iv), $\tau _0=$$\max \{n\mid
r(n)=0\}$ by (iii) and (ii), and $\tilde h(n)=h(\tau _0)$ by (iv).

(vii) For $n\geq \nu _{s+1}$ we have $t(n)=+\infty $. This means that there
exists a sequence of integers $m_i\rightarrow +\infty $ such that%
$$
\tilde h(n)=h(m_i)-(m_i-n)r(m_i)=(\rho m_i+\sigma )-(m_i-n)\rho =\rho
n+\sigma .
$$

\subsection{Filtering FDH Functions.\label{filter}}

Let $h$ be an FDH function. We will say that $h$ is {\em filterable at }$m$
if the associated function $r$ satisfies $r(n)\geq r(m)$ for all $n\geq m$.
A {\em filtration} of $h$ is a sequence of integers $m_0<m_1<\cdots <m_s$ at
which $h$ is filterable. Given such a filtration we decompose $h$ into a sum
of $s+2$ function $h_0,\ldots ,h_{s+1}$ defined as follows. We set $%
m_{-1}=-\infty $ and $m_{s+1}=+\infty $. Then the second difference $\Delta
^2h(n)=a(n)-b(n)$ may be decomposed by%
$$
a_i(n)-b_i(n):=\left\{
\begin{array}{ll}
a(n)-b(n) & \text{if }m_{i-1}<n\leq m_i, \\ 0 & \text{otherwise,}
\end{array}
\right.
$$
with $r_i(n)=\sum_{m\leq n}\left\{ a_i(m)-b_i(m)\right\} $ and $%
h_i(n)=\sum_{m\leq n}r_i(n)$ defined as in (\ref{r(n)}) and (\ref{h(n)}). If
we write $H_i(n):=h(m_i)+(n-m_i)r(m_i)$, then the $r_i(n)$ and $h_i(n)$
satisfy%
\begin{eqnarray*}
& r_i(n)  =  \left\{
\begin{array}{ll}
0 & \text{if }n\leq m_{i-1} \\ r(n)-r(m_{i-1}) & \text{if }m_{i-1}<n\leq
m_i, \\ r(m_i)-r(m_{i-1}) & \text{if }n>m_i.
\end{array}
\right. \label{ri(n)} & \\
& h_i(n)  =  \left\{
\begin{array}{ll}
0 & \text{if }n\leq m_{i-1} \\ h(n)-H_{i-1}(n) &
\text{if }m_{i-1}<n\leq m_i, \\ H_i(n)-H_{i-1}(n) & \text{if }n>m_i.
\end{array}
\right. \label{h_i(n)} &
\end{eqnarray*}The filterability of $h$ at the $m_i$ implies that $%
r_i(n)\geq 0$ for all $n$ and $i$. So the $h_i(n)$ are all FDH functions.

We call the functions $h_i(n)$ the {\em graded pieces} of the filtration. We
will say that a filtration is {\em trivial} if all but one of its graded
pieces vanish.

Now let us consider the associated function $t$ of (\ref{t(n)}). By Lemma
\ref{combin}(iii) the sequence $\tau _0<\tau _1<\cdots <\tau _s$ of all
distinct finite values of $t$ form a filtration of $h$ which we call the
{\em WHN filtration} (or weak Harder-Narasimhan filtration). Some of the
properties of this filtration are

\begin{lemma}
\label{WHN}Let $h$ be an FDH function, let $\tau _0<\tau _1<\cdots <\tau _s$
be the WHN filtration of $h$, and let $h_0,h_1,\ldots ,h_{s+1}$ be the
graded pieces of the filtration. For each $i$ let $\nu _i=\min \{n\mid
t(n)=\tau _i\}$. Then

(i)\quad The FDH function $h_0$ is torsion. It vanishes if $h$ is
torsion-free.

(ii)\quad The FDH functions $h_1,h_2,\ldots ,h_{s+1}$ are torsion-free.

(iii)\quad For $i=1,\ldots ,s$ the function $t_i$ associated to $h_i$ by (%
\ref{t(n)}) satisfies $t_i(n)=\tau _{i-1}$ for $n<\nu _i$, and $t(n)=+\infty
$ for $n\geq \nu _i$. Thus the $h_i$ are FDH functions with trivial WHN
filtrations.
\end{lemma}

\TeXButton{Proof}{\proof}(i) This is a direct translation of Lemma \ref
{combin}(vi).

(ii) This follows directly from Lemma \ref{combin}(ii) and the formula for $%
r_i(n)$.

(iii) We first suppose $n<\nu _i$. We will compute $\widetilde{h_i}(n)$ and $%
t_i(n)$ according to the definitions (\ref{htilde(n)}) and (\ref{t(n)}).
This means first computing $h_i(m)+(n-m)r_i(m)$ for all $m\geq n$.

If $m\leq \tau _{i-1}$, then $h_i(m)=r_i(m)=0$ and so $h_i(m)+(n-m)r_i(m)=0$.

If $\tau _{i-1}<m\leq \tau _i$, then%
$$
h_i(m)+(n-m)r_i(m)=\left\{ h(m)+(n-m)r(m)\right\} -\left\{ h(\tau
_{i-1})+(n-\tau _{i-1})r(\tau _{i-1})\right\} .
$$
But since $m>\tau _{i-1}$, the definitions of $\tilde h(n)$ and $t(n)$ imply
that the right side of this equation is negative for all $n$ such that $%
t(n)=\tau _{i-1}$, including $n=\nu _i-1$. The right hand side is also
linear in $n$ with slope $r(m)-r(\tau _{i-1})$ which is positive by Lemma
\ref{combin}(ii), so it must be negative for all $n<\nu _i$. Thus $%
h_i(m)+(n-m)r_i(m)<0$ if $\tau _{i-1}<m\leq \tau _i$.

If $m\geq \tau _i$, then $h_i(m)+(n-m)r_i(m)=h_i(\tau _i)+(n-\tau
_i)r_i(\tau _i)<0$ because $h_i$ is linear in this range.

So by the definitions we have $\widetilde{h_i}(n)=0$ and $\nu _i(n)=\tau
_{i-1}$ for $n<\nu _i$.

Now we suppose that $\nu _i\leq n\leq \tau _i$. Then after subtracting $%
h(\tau _{i-1})+(n-\tau _{i-1})r(\tau _{i-1})$ from both sides of the
inequality of Lemma \ref{combin}(v), we see that for all $n\leq m\leq \tau
_i $ we have
\begin{equation}
\label{tau(i)}h_i(\tau _i)+(n-\tau _i)r_i(\tau _i)\geq h_i(m)+(n-m)r_i(m).
\end{equation}
And for all $m\geq \tau _i$ we have equality in (\ref{tau(i)}) because $h_i$
is linear in this range. So by the definitions, $\tilde
h_i(n)=h_i(m)+(n-m)r_i(m)$ for all $m\geq \tau _i$, and $t_i(n)=+\infty $.

Finally if $n\geq \tau _i$, then for all $m\geq n$ we have equality in (\ref
{tau(i)}), so again we have $t_i(n)=+\infty $. \TeXButton{qed}{\qed \medskip}

\subsection{Torsion-free FDH functions $h$ with trivial WHN-filtrations.}

We wish to decompose an $h$ of this type in a certain way. For $n\gg 0$ the
function $h(n)$ is linear, so we may write it in the form $\rho (n-\nu
)+\beta $ with $\rho $, $\nu $, and $\beta $ integers such that $0\leq \beta
<\rho $. But then if $n<\nu _1$ we have by Lemma \ref{combin}(vi) that $%
h(n)=\tilde h(n)=0$ and $t(n)=\max \{n\mid h(n)=0\}$. If $n\geq \nu _1$ then
by Lemma \ref{combin}(vii) we have $t(n)=+\infty $ and $0\leq h(n)\leq $$%
\tilde h(n)=\rho (n-\nu )+\beta $. Moreover, $\nu _1=\nu $ by Lemma \ref
{combin}(iv).

We now define%
\begin{eqnarray}
\gamma^i(n) & = & \left\{
\begin{array}{ll}
(n-\nu +1)_{+} & \text{for }i=1,\ldots ,\beta , \\
(n-\nu )_{+} & \text{for }i=\beta +1,\ldots ,\rho ,
\end{array}
\right.  \\
h^i(n) & = & \min \left\{ \gamma^i(n),\left[ h(n)-\sum_{k=1}^{i-1}\gamma
^k(n)\right] _{+}\right\}. \label{hi(n)}
\end{eqnarray}These functions have the following properties:

\begin{lemma}
\label{hi}Suppose $h$ is a torsion-free FDH function such that the
associated function $t$ of (\ref{t(n)}) takes only two distinct values. For $%
i=1,\ldots ,\rho $ let $h^i{:}~{\bf Z}\rightarrow {\bf Z}_{\geq 0}$ be the
function defined in (\ref{hi(n)}). Then

(i)\quad The $h^i$ satisfy $\sum_{i=1}^\rho h^i=h$,

(ii)\quad The $h^i$ are torsion-free FDH functions of rank $1$. The degree
of $h^i$ is $-\nu $ for $i=1,\ldots ,\beta $, and $-\nu -1$ for $i=\beta
+1,\ldots ,\rho $.

(iii)\quad The deficiency of $h^i$ is positive for $i=1,\ldots ,\beta $.
\end{lemma}

\TeXButton{Proof}{\proof}(i) First note that if $n<\nu $, then $\gamma
^i(n)=h^i(n)=0$ for all $i$. But $h(n)=0$ as well. So this case is fine. If $%
n\geq \nu $, then we have $0\leq h(n)\leq \tilde h(n)=\sum_{i=1}^\rho \gamma
^i(n)$. So there exists a $k$ such that $\sum_{i=1}^{k-1}\gamma ^i(n)\leq
h(n)\leq \sum_{i=1}^k\gamma ^i(n)$. Then

\begin{equation}
\label{h(i)(n)}h^i(n)=\left\{
\begin{array}{ll}
\gamma ^i(n) & \text{for }i=1,\ldots ,k-1 \\ h(n)-\sum_{i=1}^{k-1}\gamma
^i(n) & \text{for }i=k \\ 0 & \text{for }i=k+1,\ldots ,\rho
\end{array}
\right.
\end{equation}
and the sum is $h(n)$. This completes the proof of (i).

For (ii) we first show that the $h^i$ are torsion-free FDH functions, i.e.\
that $\Delta h^i(n)>0$ for all $n$ such that $h^i(n)>0$. To verify this, we
may clearly assume that $n\geq \nu $ since otherwise $h^i(n)=0$. Now note
that if one has a function of the form $h^i=\min (f,g)$, then in order to
show that $h^i(n)>0$ implies $\Delta h^i(n)>0$ it is enough to show that $%
f(n)>0$ implies $\Delta f(n)>0$ and that $g(n)>0$ implies $\Delta g(n)>0$.
So now consider the case $i=1,\ldots ,\beta $. The function $f(n):=\gamma
^i(n)$ satisfies $\Delta \gamma ^i(n)=1>0$. And if $g(n):=\left[
h(n)-\sum_{k=1}^{i-1}\gamma ^k(n)\right] _{+}>0$, then $h(n)>(i-1)(n-\nu +1)$%
. But from the definition of $\tilde h(n)$ we have%
$$
0=\tilde h(\nu -1)\geq h(n)-(n-\nu +1)\Delta h(n).
$$
So $(n-\nu +1)\Delta h(n)\geq h(n)$. Thus $\Delta h(n)>i-1$, and $\Delta
g(n)>0$. This proves that $h^i$ is a torsion-free FDH function for $%
i=1,\ldots ,\beta $.

The proof that $h^i$ is a torsion-free FDH function for $i=\beta +1,\ldots
,\rho $ is similar except that one uses $\tilde h(\nu )=\beta $ to obtain $%
(n-\nu +1)\Delta h(n)\geq h(n)-\beta $.

For the rank and degree of $h^i$ note that for $n\gg 0$ the formula (\ref
{hi(n)}) becomes%
$$
h^i(n)=\min \left\{ \gamma ^i(n),\sum_{k=i}^\rho \gamma ^k(n)\right\}
=\gamma ^i(n)=\left\{
\begin{array}{ll}
n-\nu +1 & \text{for }i=1,\ldots ,\beta , \\ n-\nu & \text{for }i=\beta
+1,\ldots ,\rho .
\end{array}
\right.
$$

For (iii) note that for $i=1,\ldots ,\beta $ the deficiency of $h^i$ is $%
\sum_{n\geq \nu }(n-\nu +1-h^i(n))$. All the terms in this sum are
nonnegative, so it is enough to show that $h^i(\nu )=0$. But recall that $%
t(\nu -1)=\max \{n\mid h(n)=0\}$. And by Lemma \ref{combin}(i) $t(\nu
-1)\geq \nu $. Thus $h(\nu )=0$, from which $h^i(\nu )=0$ for all $i$ by
(i). Part (iii) now follows. \TeXButton{qed}{\qed \medskip}

As a final numerical result we wish to compare the decompositions of Lemmas
\ref{WHN} and \ref{hi}. To do this we introduce an order on torsion-free FDH
functions. Namely if $h$ and $h^{\prime }$ are torsion-free FDH functions of
ranks $\rho $ and $\rho ^{\prime }$, degrees $d$ and $d^{\prime }$, and
deficiencies $\delta $ and $\delta ^{\prime }$, then $h\succeq h^{\prime }$
(resp.\ $h\succ h^{\prime }$) if $d/\rho >d^{\prime }/\rho ^{\prime }$ or if
$d/\rho =d^{\prime }/\rho ^{\prime }$ and $\delta /\rho \leq \delta ^{\prime
}/\rho ^{\prime }$ (resp.\ $\delta /\rho <\delta ^{\prime }/\rho ^{\prime }$%
).

\begin{lemma}
\label{compar}Let $h$ be an FDH function, let $h_0,h_1,\ldots ,h_{s+1}$ be
the graded pieces of the WHN filtration of $h$ of Lemma \ref{WHN}. For $%
i=1,\ldots ,s+1$, let $\rho _i$ denote the rank of $h_i$, and let $%
h_i=\sum_{j=1}^{\rho _i}h_i^j$ be the decomposition of $h_i$ of Lemma \ref
{hi}. Then

(i)\quad $h_i^1\succeq h_i^2\succeq \cdots \succeq h_i^{\rho _i}$ for $%
i=1,\ldots ,s+1$, and

(ii)\quad $h_i^{\rho _i}$ $\succ h_{i+1}^1$ for $i=1,\ldots ,s$.
\end{lemma}

\TeXButton{Proof}{\proof}First we introduce some notation. As in Lemma \ref
{hi} if $n\gg 0$, then we may write $h_i(n)=\rho _i(n-\nu _i)+\beta _i$ with
$0\leq \beta _i<\rho _i$. For each $i$ and $j=1,\ldots ,\rho _i$ we define%
$$
\begin{array}{lclcl}
\eta _i^j & := & \min \{n\mid h_i^j(n)>0\} & = & \min \{n\mid
h_i(n)>\sum_{k=1}^{j-1}\gamma _i^k(n)\}, \\
\zeta _i^j & := & \min \{n\mid h_i^j(n)=\gamma _i^j(n)>0\} & = & \min
\{n\mid h_i(n)\geq \sum_{k=1}^j\gamma _i^k(n)>0\}.
\end{array}
$$
where $\gamma _i^k$ is as in Lemma \ref{hi}. Then if $1\leq j\leq \beta _i$
(resp.\ if $\beta _i+1\leq j\leq \rho _i$), the FDH function $h_i$ has rank $%
1$, degree $d_i^j=-\nu _i$ (resp.\ $d_i^j=-\nu _i-1$), and deficiency $%
\delta _i^j$ satisfying $\binom{\eta _i^j+d_i^j+1}2\leq \delta _i^j\leq
\binom{\zeta _i^j+d_i^j+1}2$.

(i) For all $i$ and $j=1,\ldots ,\rho _i-1$ we have $\zeta _i^j\leq \eta
_i^{j+1}$. Hence $h_i^1\succeq h_i^2\succeq \cdots \succeq h_i^{\beta _i}$
for all $i$ because all these functions have the same rank, the same degree $%
\alpha _i$, and nondecreasing deficiencies. Similarly $h_i^{\beta
_i+1}\succeq h_i^{\beta _i+2}\succeq \cdots \succeq h_i^{\beta _i}$. And $%
h_i^{\beta _i}\succ h_i^{\beta _i+1}$ by reason of degree.

(ii) Note that by Lemmas \ref{combin} (v) and \ref{WHN}, the FDH function $%
h_i^{\rho _i}$ has degree $-\nu _i-1$ and has
$$
\zeta _i^{\rho _i}=\min \{n>\nu _i\mid h(n)=h(\tau _i)+(n-\tau _i)r(\tau
_i)\}\leq \tau _i-1,
$$
while $h_{i+1}^1$ has degree $-\nu _{i+1}$ or $-\nu _{i+1}-1$ and has%
$$
\eta _{i+1}^1=\min \{n>\nu _i\mid h(n)>h(\tau _i)+(n-\tau _i)r(\tau
_i)\}=\tau _i+1.
$$
Since $\nu _i<\nu _{i+1}$, the degree of $h_i^{\rho _i}$ is at least that of
$h_{i+1}^1$, and in case of equality the former function has a smaller
deficiency than the latter.\TeXButton{qed}{\qed \medskip}

\subsection{Effective computation of the $\tau _i$ and $\nu _i$.\label
{effective}}

The $\tau _i$ and $\nu _i$ defined in Lemma \ref{combin} and referred to in
the statement of Theorem \ref{main} may be effectively computed from $h(n)$
or $r(n)=\Delta h(n)$. Note that $\Delta r(n)$ is $0$ for all but finitely
many $n$. For $t=+\infty $ we write $r(t)=\rho $ and $h(t)-tr(t)=\sigma $.

According to Lemma \ref{combin}(ii),
$$
\{\tau _i\}_{i=0}^{s+1}\subset T:=\left\{ n\mid r(m)>r(n)\text{ for all }%
m>n\right\} \cup \{+\infty \}.
$$
The set $T$ may be computed by passing through the finite set%
$$
T^{\prime }:=\left\{ n\mid \Delta r(n+1)>0\right\} \cup \{+\infty \}
$$
in descending order and purging those $n\in T^{\prime }$ such that $r(n)\geq
r(m)$ where $m$ is the smallest unpurged element of $T^{\prime }$ larger
than $n$. The minimal element of $T$ is $\tau _0$ since it is the unique $%
t\in T$ such that $r(t)=0$. The other $\tau _i$ and $\nu _i$ may be computed
recursively as follows.

Suppose we have computed $\tau _0,\ldots ,\tau _{i-1}$ and $\nu _1$,$\ldots
,\nu _{i-1}$. We now need to find for which $x>\nu _{i-1}$ there is a $t\in
T $ with $t>\tau _{i-1}$ such that%
$$
h(t)+(x-t)r(t)\geq h(\tau _{i-1})+(x-\tau _{i-1})r(\tau _{i-1}).
$$
For $t=+\infty $, the left side should be read as $\rho x+\sigma $ according
to our above conventions. Each of the inequalities is equivalent to%
$$
x\geq x_t:=\frac{\left( h(t)-tr(t)\right) -\left( h(\tau _{i-1})-\tau
_{i-1}r(\tau _{i-1})\right) }{r(t)-r(\tau _{i-1})}
$$
So $\nu _i=\min \left\{ \left\lceil x_t\right\rceil \mid t\in T\text{ and }%
t>\tau _{i-1}\right\} $ where the notation $\left\lceil x_t\right\rceil $
means the smallest integer greater than or equal to $x_t$. We then look at
those $t$ such that $\left\lceil x_t\right\rceil =\nu _i$, and pick out
those among them for which $h(t)+(\nu _i-t)r(t)$ is maximal. The largest of
these $t$ is $\tau _i$.

We continue until some $\tau _i=+\infty $.

The other invariants in the statement of Theorem \ref{main} and the proof of
Lemma \ref{compar} may be computed as%
\begin{eqnarray*}
\rho _i & = & r(\tau _i)-r(\tau _{i-1}), \\
\beta _i & = & \left( h(\tau _i)+(\nu _i-\tau _i)r(\tau _i)\right) -\left(
h(\tau _{i-1})+(\nu _i-\tau _{i-1})r(\tau _{i-1})\right) .
\end{eqnarray*}

\section{A Special Filtration on ${\bf P}^2$\label{sec2}}

In this section we proceed to give sheaf-theoretic significance to the
numerical computations of the previous section. We do this by introducing
the WHN filtration on the general torsion-free sheaf ${\cal E}$ with a given
Hilbert function. The Hilbert function of the graded pieces ${\rm gr}_i(%
{\cal E})$ of the filtration are those given by Lemma \ref{WHN}. Lemma \ref
{hi} is then used to show that the graded pieces satisfy ${\rm Hom}({\rm gr}%
_i({\cal E}),{\rm gr}_j({\cal E})(-1))=0$ for all $i\leq j$. Lemma \ref
{compar} is used to show that the WHN filtration is compatible with the
Harder-Narasimhan filtration for Gieseker-Maruyama stability.

\subsection{Hilbert Functions.}

Recall that any coherent sheaf ${\cal E}$ on ${\bf P}^2$ without
zero-dimensional torsion has a free resolution of the form
\begin{equation}
\label{phi}0\rightarrow \bigoplus_n{\cal O}_{{\bf P}^2}(-n)^{b(n)}\stackrel{%
\phi }{\rightarrow }\bigoplus_n{\cal O}_{{\bf P}^2}(-n)^{a(n)}\rightarrow
{\cal E}\rightarrow 0.
\end{equation}
The $a(n)$ and $b(n)$ are related to the Hilbert function $n\mapsto h^0(%
{\cal E}(n))$ of ${\cal E}$ via
\begin{equation}
\label{Hilbert}\sum_nh^0({\cal E}(n))\,t^n=(1-t)^{-3}\sum_n\{a(n)-b(n)\}t^n
\end{equation}
or $a(n)-b(n)=\Delta ^3h^0({\cal E}(n))$. So the $a(n)$ and $b(n)$ determine
the Hilbert function of ${\cal E}$ which conversely determines the
differences $a(n)-b(n)$. If ${\cal E}$ is sufficiently general, then for all
$n$ either $a(n)=0$ or $b(n)=0$ according to the sign of $a(n)-b(n)$, so the
Hilbert function then actually determines the $a(n)$ and $b(n)$.

The next theorem is the filtered Bertini theorem as applied to the special
case of ${\bf P}^2$:

\begin{theorem}
{\rm \label{Chang}(Chang \cite{Ch})} A general map $\phi {:}\ \bigoplus_n%
{\cal O}_{{\bf P}^2}(-n)^{b(n)}\rightarrow \bigoplus_n{\cal O}_{{\bf P}%
^2}(-n)^{a(n)}$ is injective if and only if the function $h$ whose
Poincar\'e series is%
$$
\sum_nh(n)t^n=(1-t)^{-2}\sum_n\{a(n)-b(n)\}t^n
$$
is an FDH function. Moreover, the cokernel ${\cal E}$ of a general $\phi $ is

$\bullet $ torsion-free with at worst singular points of multiplicity $1$ if
$h$ is torsion-free,

$\bullet $ locally free if $h$ is locally free,

$\bullet $ a line bundle on a curve with normal crossings if $h$ is torsion,

$\bullet $ a line bundle on a smooth curve if $h$ is torsion and
unfilterable.
\end{theorem}

Here the {\em multiplicity} of a singular point $P$ of a torsion-free sheaf $%
{\cal E}$ on a smooth surface is the length of ${\cal E}_P^{\vee \vee }/%
{\cal E}_P$.

Theorem \ref{Chang} and formula (\ref{Hilbert}) allow us to define the {\em %
FDH function of a coherent sheaf }${\cal E}$ without zero-dimensional
torsion as $h:=\Delta h^0({\cal E}(n))$. The rank of the FDH function $h$ is
then ${\rm rk}({\cal E})$, and the degree of $h$ is $c_1({\cal E)}$. If this
rank is $\rho $, and the degree is written as $\rho \alpha +\beta $ with $%
\alpha $ and $\beta $ integers such that $0\leq \beta <\rho $, then the
deficiency of $h$ is $c_2({\cal E})-c_2({\cal F})$ where ${\cal F}={\cal O}_{%
{\bf P}^2}(\alpha +1)^\beta \oplus {\cal O}_{{\bf P}^2}(\alpha )^{\rho
-\beta }$.

We may speak of a generic sheaf with FDH function $h$ because of the
following fact, which is well known and which we therefore state without
proof:

\begin{lemma}
The coherent sheaves without zero-dimensional torsion on ${\bf P}^2$ with a
fixed Hilbert function or FDH function form an irreducible and smooth
locally closed substack of the stack of coherent sheaves on ${\bf P}^2$.
\end{lemma}

Our next lemma relates the filterability of the FDH function $h$ to sheaf
theory. It is essentially a generalization of Davis' decomposition lemma
\cite{D} to higher rank.

\begin{lemma}
\label{gener}Suppose $h$ is an FDH function which is filterable at an
integer $m$. Suppose ${\cal E}$ is a general coherent sheaf without
zero-dimensional torsion with FDH function $h$. Let ${\cal E}_{\leq m}$ be
the subsheaf of ${\cal E}$ generated by $H^0({\cal E}(m))$. Then ${\cal E}%
_{\leq m}$ has FDH function $h_m^{{\rm sub}}$, and ${\cal E}/{\cal E}_{\leq
m}$ has FDH function $h_m^{{\rm quot}}$. Moreover, if ${\cal E}$ is generic,
then so are ${\cal E}_{\leq m}$ and ${\cal E}/{\cal E}_{\leq m}$. If ${\cal E%
}$ is the cokernel of an injection $\phi {:}\ \bigoplus_n{\cal O}_{{\bf P}%
^2}(-n)^{b(n)}\rightarrow \bigoplus_n{\cal O}_{{\bf P}^2}(-n)^{a(n)}$, then $%
{\cal E}_{\leq m}$ and ${\cal E}/{\cal E}_{\leq m}$ have resolutions%
\begin{eqnarray*}
& 0\rightarrow \bigoplus_{n\leq m}{\cal O}_{{\bf P}^2}(-n)^{b(n)}\rightarrow
\bigoplus_{n\leq m}{\cal O}_{{\bf P}^2}(-n)^{a(n)}\rightarrow {\cal E}_{\leq
m}\rightarrow 0 &  \\
& 0\rightarrow \bigoplus_{n>m}{\cal O}_{{\bf P}^2}(-n)^{b(n)}\rightarrow
\bigoplus_{n>m}{\cal O}_{{\bf P}^2}(-n)^{a(n)}\rightarrow {\cal E}%
_{>m}\rightarrow 0 &
\end{eqnarray*}
\end{lemma}

\TeXButton{Proof}{\proof}Consider the morphism of exact sequences%
$$
\begin{array}{ccccccl}
0\rightarrow & \bigoplus_{n\leq m}{\cal O}_{{\bf P}^2}(-n)^{b(n)} &
\rightarrow & \bigoplus {\cal O}_{{\bf P}^2}(-n)^{b(n)} & \rightarrow &
\bigoplus_{n>m}{\cal O}_{{\bf P}^2}(-n)^{b(n)} & \rightarrow 0 \\
& \downarrow \phi ^{\prime \prime } &  & \downarrow \phi &  & \downarrow
\phi ^{\prime } &  \\
0\rightarrow & \bigoplus_{n\leq m}{\cal O}_{{\bf P}^2}(-n)^{a(n)} &
\rightarrow & \bigoplus {\cal O}_{{\bf P}^2}(-n)^{a(n)} & \rightarrow &
\bigoplus_{n>m}{\cal O}_{{\bf P}^2}(-n)^{a(n)} & \rightarrow 0.
\end{array}
$$
The Poincar\'e series associated to $\phi ^{\prime }$ and $\phi ^{\prime
\prime }$ are, respectively,%
\begin{eqnarray*}
(1-t)^{-2}\sum_{n>m}\{a(n)-b(n)\}t^n & = & \sum_nh_m^{\rm quot}(n)t^n, \\
(1-t)^{-2}\sum_{n\leq m}\{a(n)-b(n)\}t^n & = & \sum_nh_m^{\rm sub}(n)t^n.
\end{eqnarray*}Since $h$ is filterable at $m$, $h_m^{{\rm quot}}$ and $h_m^{%
{\rm sub}}$ are both FDH functions by (\ref{filter}). If $\phi $ is general,
then $\phi ^{\prime }$ is general and so injective by the filtered Bertini
Theorem \ref{Chang}. In any case the snake lemma yields an exact sequence%
$$
0\rightarrow \ker (\phi ^{\prime })\rightarrow {\rm cok(\phi ^{\prime \prime
})}\stackrel{\psi }{\rightarrow }{\cal E}\rightarrow {\rm cok(\phi ^{\prime
})}\rightarrow 0
$$
such that ${\rm im}(\psi )={\cal E}_{\leq m}$. So if $\phi $ is general,
then ${\rm cok(\phi ^{\prime \prime })}={\cal E}_{\leq m}$ and ${\rm %
cok(\phi ^{\prime })}={\cal E}/{\cal E}_{\leq m}$, and they have FDH
functions $h_m^{{\rm sub}}$ and $h_m^{{\rm quot}}$, respectively. Finally,
if ${\cal E}$ is generic, then $\phi $ is generic, which implies the
genericity of $\phi ^{\prime }$ and $\phi ^{\prime \prime }$ and thus of $%
{\cal E}_{\leq m}$ and ${\cal E}/{\cal E}_{\leq m}$. \TeXButton{qed}
{\qed \medskip}

\subsection{Davis' Decomposition Lemma.}

If ${\cal E}$ is torsion-free and $r(m)=1=\min _{n\geq m}r(n)$, then the
lemma holds without any condition that ${\cal E}$ be general. This is
because the vanishing of $\ker (\phi ^{\prime })$ may be shown without
invoking the filtered Bertini theorem. For the condition $r(m)=1$ implies
that ${\rm cok(\phi ^{\prime \prime })}$ is of rank $1$. Since it has
nonzero image in ${\cal E}$ which is torsion-free, it follows that $\ker
(\phi ^{\prime })$ is of rank $0$. But since $\ker (\phi ^{\prime })\subset
\bigoplus_{n>m}{\cal O}_{{\bf P}^2}(-n)^{b(n)}$, it is also torsion-free. So
it vanishes. In the case of an ${\cal E}$ of rank one, this is more or less
Davis' Decomposition Lemma \cite{D}.

\subsection{The WHN Filtration of a General Sheaf.\label{WHNfilt}}

Because of the last lemma, if $h$ is an FDH function with a filtration $%
m_0<m_1<\cdots <m_s$, then the general sheaf ${\cal E}$ with FDH function $h$
will have a filtration
\begin{equation}
\label{filt.sheaf}0\subset {\cal E}_{\leq m_0}\subset {\cal E}_{\leq
m_1}\subset \cdots \subset {\cal E}_{\leq m_s}\subset {\cal E.}
\end{equation}
If we write ${\rm gr}_i({\cal E})={\cal E}_{\leq m_i}/{\cal E}_{\leq
m_{i-1}} $ for $i=1,\ldots ,s$, and ${\rm gr}_0({\cal E})={\cal E}_{\leq
m_0} $ and ${\rm gr}_{s+1}({\cal E})={\cal E}/{\cal E}_{\leq m_s}$, then the
FDH function of ${\rm gr}_i({\cal E})$ is the function $h_i$ of (\ref{filter}%
). If ${\cal E}$ has resolution%
$$
0\rightarrow \bigoplus_n{\cal O}_{{\bf P}^2}(-n)^{b(n)}\rightarrow
\bigoplus_n{\cal O}_{{\bf P}^2}(-n)^{a(n)}\rightarrow {\cal E}\rightarrow 0,
$$
then each ${\rm gr}_i({\cal E})$ has resolution%
$$
0\rightarrow \bigoplus_{m_{i-1}<n\leq m_i}{\cal O}_{{\bf P}%
^2}(-n)^{b(n)}\rightarrow \bigoplus_{m_{i-1}<n\leq m_i}{\cal O}_{{\bf P}%
^2}(-n)^{a(n)}\rightarrow {\rm gr}_i({\cal E})\rightarrow 0.
$$

If we apply this with the WHN filtration of $h$ of Lemma \ref{WHN}, then we
call the resulting filtration of ${\cal E}$ the {\em WHN filtration }of the
sheaf ${\cal E}$. This filtration only exists for a general sheaf with FDH
function $h$ because the construction of the filtration of the sheaf
depended ultimately on the filtered Bertini theorem.\bigskip\

We now recall some terminology. If ${\cal E}$ is a coherent sheaf of rank $%
\rho >0$ on ${\bf P}^2$, then its reduced Hilbert polynomial is $P_{{\cal E}%
}(n):=\chi ({\cal E}(n))/\rho $. Such polynomials may be ordered by $P_{%
{\cal E}}\succ P_{{\cal F}}$ (resp.$\ P_{{\cal E}}\succeq P_{{\cal F}}$) if $%
P_{{\cal E}}(n)>P_{{\cal F}}(n)$ (resp.\ $P_{{\cal E}}(n)\geq P_{{\cal F}%
}(n) $) for $n\gg 0$. This order is compatible with the order on FDH
functions of Lemma \ref{compar} in the sense that if ${\cal E}$ has FDH
function $h_{{\cal E}}$ and ${\cal F}$ has FDH function $h_{{\cal F}}$, then
$P_{{\cal E}}\succ P_{{\cal F}}$ if and only if $h_{{\cal E}}\succ h_{{\cal F%
}}$, and $P_{{\cal E}}\succeq P_{{\cal F}}$ if and only if $h_{{\cal E}%
}\succeq h_{{\cal F}}$.

\begin{lemma}
\label{ideals}Let $h$ be an FDH function, let the $h_i$, $\rho _i$, $\nu _i$%
, $\beta _i$, and $h_i^j$ be as in Lemma \ref{compar}. Then there exists a
coherent sheaf ${\cal F}$ with FDH function $h$ of the form ${\cal F}%
=\bigoplus_{i=0}^{s+1}{\cal F}_i$ such that

(i)$\quad {\cal F}$ admits the WHN filtration with graded pieces ${\rm gr}_i(%
{\cal F})\cong {\cal F}_i$.

(ii)\quad ${\cal F}_0$ is the torsion subsheaf of ${\cal F}$.

(iii)\quad For $i=1,\ldots ,s+1$ we have ${\cal F}_i=\bigoplus_{j=1}^{\beta
_i}{\cal I}_{Z_i^j}(-\nu _i)\oplus \bigoplus_{j=\beta _i+1}^{\rho _i}{\cal I}%
_{Z_i^j}(-\nu _i-1)$ where the $Z_i^j$ are disjoint sets of distinct points.
Moreover, $Z_i^j\neq \emptyset $ for $1\leq j\leq \beta _i$.
\end{lemma}

\TeXButton{Proof}{\proof}Let ${\cal F}_0$ be a general sheaf with FDH
function $h_0$, and let ${\cal F}_i^j$ be a general sheaf with FDH function $%
h_i^j$. For $i=1,\ldots ,s+1$, let ${\cal F}_i=\bigoplus_{j=1}^{\rho _i}%
{\cal F}_i^j$ and ${\cal F}=\bigoplus_{i=0}^{s+1}{\cal F}_i$. Then ${\cal F}$
has FDH function $h_0+\sum_{i=1}^{s+1}\sum_{j=1}^{\rho _i}h_i^j=h$. We now
verify the three asserted conditions in reverse order.

(iii) By Lemma \ref{hi} (ii), the $h_i^j$ are torsion-free sheaves of rank $%
1 $. So by Theorem \ref{Chang} ${\cal F}_i^j$ is a twist of an ideal sheaf
of a set $Z_i^j$ of distinct points. The twist is given by the degree of $%
h_i^j$, which is $-\nu _i$ if $1\leq j\leq \beta _i$ (resp.\ $-\nu _i-1$ if $%
\beta _i+1\leq j\leq \rho _i$). Replacing the $Z_i^j$ by projectively
equivalent sets of points if necessary, we may assume that the $Z_i^j$ are
disjoint. Finally, the cardinality of $Z_i^j$ is the deficiency of $h_i^j$,
which is positive if $1\leq j\leq \beta _i$ by Lemma \ref{hi} (iii).

(ii) The sheaf ${\cal F}_0$ is torsion because the function $h_0$ is torsion
by Lemma \ref{WHN} (i). The other factors ${\cal F}_i$ in the direct sum $%
{\cal F}$ are torsion-free by part (iii) which we just proved. So ${\cal F}%
_0 $ is exactly the torsion subsheaf of ${\cal F}$.

(i) We need to show that for $i=0,\ldots ,s$, the subsheaf ${\cal F}_{\leq
\tau _i}\subset {\cal F}$ is $\bigoplus_{k=0}^i{\cal F}_k$. So we show that $%
{\cal F}_{i,\leq \tau _i}={\cal F}_i$ and ${\cal F}_{i+1,\leq \tau _i}=0$.

First suppose that $g$ is an FDH function of rank $\rho $ and degree $d$,
and if $\tau $ is an integer such that $g(n)=\rho (n+1)+d$ for all $n\geq
\tau -1$, then $g$ is filterable at $\tau $, and $g_\tau ^{{\rm sub}}=g$ and
$g_\tau ^{{\rm quot}}=0$. So if ${\cal G}$ is a general sheaf with FDH
function $g$, then ${\cal G}_{\leq \tau }={\cal G}$. If we apply this with $%
g=h_0$ and $\tau =\tau _0$, we see that ${\cal F}_{0,\leq \tau _0}={\cal F}%
_0 $. We may also apply it with $g=h_i^j$ and $\tau =\tau _i$ because $\tau
_i\geq \zeta _i^j+1$ where $\zeta _i^j$ is as in the proof of Lemma \ref
{compar}. Thus for $i=1,\ldots ,s$, we have ${\cal F}_{i,\leq \tau
_i}=\bigoplus_{j=1}^{\rho _i}{\cal F}_{i,\leq \tau _i}^j=\bigoplus {\cal F}%
_i^j={\cal F}_i$.

To show that ${\cal F}_{i+1,\leq \tau _i}=0$ we need to show that all $%
h_{i+1}^j(\tau _i)=0$. But if $\eta _{i+1}^j$ is as in the proof of Lemma
\ref{compar}, then $h_{i+1}^j(n)$ for all $n<\eta _{i+1}^j$. But as we
showed there $\eta _{i+1}^j\geq \eta _{i+1}^1\geq \tau _i+1$. This completes
the proof of the lemma.\TeXButton{qed}{\qed \medskip}

\paragraph{Proof of Theorem \ref{main}.}

Part (i) was shown in (\ref{WHNfilt}). Parts (ii), (iv), and (v) then
describe open properties, so it is enough to verify them for the sheaf $%
{\cal F}$ of Lemma \ref{ideals}. Part (ii) then follows from Lemma \ref
{ideals}(ii). To derive (iv) from \ref{ideals}(iii), note that the fact that
the $Z_i^k$ are all disjoint implies that ${\cal H}om({\cal I}_{Z_i^k},{\cal %
I}_{Z_j^l})\cong {\cal I}_{Z_j^l}$. So ${\rm Hom}({\rm gr}_i({\cal E}),{\rm %
gr}_j({\cal E})(-1))$ is a sum of terms of the forms $H^0({\cal I}%
_{Z_j^l}(\nu _i-\nu _j))$, $H^0({\cal I}_{Z_j^l}(\nu _i-\nu _j-1))$ and $H^0(%
{\cal I}_{Z_j^l}(\nu _i+1-\nu _j))$. In the first two forms the  cohomology
vanishes because the twists $\nu _i-\nu _j$ or $\nu _i-\nu _j-1$ are
negative. For the third form the twist $\nu _i+1-\nu _j$ is nonpositive, but
even if it is zero $H^0({\cal I}_{Z_j^l})$ vanishes because this form only
occurs with $1\leq l\leq \beta _j$ and in that case $Z_j^l\neq \emptyset $.

Before beginning on (v) note that since the ${\rm gr}_i({\cal F}%
)=\bigoplus_{j=1}^{\rho _i}{\cal F}_i^j$ is a direct sum of semistable
sheaves, the graded pieces of the Harder-Narasimhan of ${\rm gr}_i({\cal F})$
are direct sums of ${\cal F}_i^j$'s with proportional Hilbert polynomials.
Hence any non-torsion quotient sheaf ${\cal G}$ of ${\rm gr}_i({\cal F})$
has $P_{{\cal G}}\succeq \min _j\{P_{{\cal F}_i^j}\}$, which is $P_{{\cal F}%
_i^{\beta _i}}$ by Lemma \ref{compar}(i), and any nonzero subsheaf ${\cal H}$
has $P_{{\cal H}}\preceq $ $\max _j\{P_{{\cal F}_i^j}\}=P_{{\cal F}_i^1}$.

Now to show that the Harder-Narasimhan filtration of ${\cal F}$ is a
refinement of the WHN filtration, we need to show that for $1\leq i\leq s$,
if ${\cal G}$ is a nonzero torsion-free quotient of ${\rm gr}_i({\cal F})$
and ${\cal H}$ a nonzero subsheaf of ${\rm gr}_{i+1}({\cal F})$, then $P_{%
{\cal G}}\succ P_{{\cal H}}$. But by the previous paragraph and Lemma \ref
{compar}(ii) we have $P_{{\cal G}}\succeq P_{{\cal F}_i^{\beta _i}}\succ P_{%
{\cal F}_{i+1}^1}\succeq P_{{\cal H}}$.

To show the second assertion of (v) we now need to show that every nonzero
subsheaf of ${\rm gr}_i({\cal F})$ has slope at most $-\nu _i$ and every
non-torsion quotient sheaf has slope at least $-\nu _i-1$. But this is now
clear.

In part (iii) the isomorphisms ${\rm gr}_i({\cal F}){\mid }_L\cong {\cal O}%
_L(-\nu _i)^{\beta _i}\oplus {\cal O}_L(-\nu _i-1)^{\rho _i-\beta _i}$
follow from Lemma \ref{ideals}(iii). Because these latter sheaves are rigid
(i.e.\ generic in the stack of coherent sheaves on $L$), a general ${\cal E}$
must have ${\rm gr}_i({\cal E}){\mid }_L$ isomorphic to ${\cal O}_L(-\nu
_i)^{\beta _i}\oplus {\cal O}_L(-\nu _i-1)^{\rho _i-\beta _i}$. We now claim
that any filtered sheaf ${\cal H}$ such that ${\rm Ext}^1({\rm gr}_i({\cal H}%
),{\rm gr}_j({\cal H}))=0$ for all $i>j$ has ${\cal H}\cong \bigoplus_i{\rm %
gr}_i({\cal H})$. This claim can easily be verified by induction on the
length of the filtration. To apply this to ${\cal E}{\mid }_L$, we need to
verify that if $i>j$, then ${\rm Ext}^1({\rm gr}_i({\cal E}){\mid }_L,{\rm gr%
}_j({\cal E}){\mid }_L)=0$. But ${\rm Ext}^1({\rm gr}_i({\cal E}){\mid }_L,%
{\rm gr}_j({\cal E}){\mid }_L)$ is a direct sum of terms of the form $H^1(%
{\cal O}_L(\nu _i-\nu _j+\epsilon ))$ with $\epsilon \in \{-1,0,1\}$. Since
the $\nu _i$ form a strictly increasing sequence of integers, the twists $%
\nu _i-\nu _j+\epsilon $ are all nonnegative$.$ So the $H^1$ vanish.
Therefore ${\cal E}{\mid }_L\cong \bigoplus_i{\rm gr}_i({\cal E}){\mid }_L$,
completing the proof of (iii). \TeXButton{qed}{\qed \medskip}

\end{document}